\documentclass[12pt]{article}

\usepackage{amssymb,amsmath}
\usepackage[margin=0.75in]{geometry}
\usepackage{mathrsfs}
\usepackage[all]{xy}
\usepackage[usenames,dvipsnames,svgnames,table]{xcolor}
\usepackage{tikz}
\usepackage{tkz-graph}
\usepackage{tkz-berge}
\usetikzlibrary{arrows,calc,automata,shadows,backgrounds,positioning,intersections,fadings,decorations.pathreplacing,shapes,snakes}
\usetikzlibrary{calc}
\usetikzlibrary{decorations.markings}
\usepackage{tikz-cd}
\usepackage{gensymb}
\usepackage[shortlabels]{enumitem}
\usepackage{float}
\usepackage{subcaption}

\usepackage[numbers]{natbib}

\newtheorem{theorem}{Theorem}[subsection]

\newtheorem{definition}[theorem]{Definition}
\newtheorem{example}[theorem]{Example}

\newtheorem{lemma}[theorem]{Lemma}

\numberwithin{equation}{section}

\def\qed{\hfill {\hbox{${\vcenter{\vbox{             
   \hrule height 0.4pt\hbox{\vrule width 0.4pt height 6pt
   \kern5pt\vrule width 0.4pt}\hrule height 0.4pt}}}$}}}

\newenvironment{proof}[1][Proof]{\smallskip\noindent{\bf #1.}\quad}
{\qed\par\medskip}

\begin{document}

\begin{center}
\Large Lewis Carroll and the Red Hot Potato:\\ a graph theoretic approach to a linear algebraic identity
\end{center}

\normalsize
\begin{center}
    Melanie Fraser
\end{center}
\textbf{Abstract.} The Lewis Carroll Identity expresses the determinant of a matrix in terms of subdeterminants obtained by deleting one row and column or a pair of rows and columns. Using the matrix tree theorem, we can convert this into an equivalent identity involving sums over pairs of forests. Unlike the Lewis Carroll Identity, the Forest Identity involves no minus signs. In 2011, Vlasev and Yeats \cite{vlasev} suggested that such a Forest Identity could be proven using edge transfers similar to Zeilberger's 1997 matrix proof. However, until now, such an algorithm has not yet been developed. In this paper, we provide this edge transfer algorithm and a bijective proof for both the Lewis Carroll Identity and Forest Identity. This bijection is implemented by the Red Hot Potato algorithm, so called because the way edges get tossed back and forth between the two forests is reminiscent of the children's game of hot potato.

\section{Introduction}

This paper presents a graph theoretic interpretation of the Lewis Carroll Identity (also known as Dodgson's rule or the Desnanot-Jacobi Identity). Most proofs of the Lewis Carroll Identity are algebraic, with the first complete proof for any $n$ given by Jacobi in 1833 \cite{muirvol1}. Zeilberger gave a combinatorial proof of this identity in 1997 \cite{zeilberger}. The Lewis Carroll Identity is largely used in Dodgson condensation, an iterative technique for evaluating determinants. More recently, the identity has been used as a basis for Dodgson polynomials, spanning forest polynomials that have been used to study Feynman graphs and perturbative quantum field theory \cite{brown}, \cite{schnetz}. In 2011, Vlasev and Yeats \cite{vlasev} suggested that a Forest Identity based on the Lewis Carroll Identity could be proven using edge transfers similar to Zeilberger's 1997 matrix proof. However, until now, such an algorithm has not yet been developed. In this paper, we provide this edge transfer algorithm and proof for both the Lewis Carroll Identity and Forest Identity.

\begin{definition}
Let $U$ and $W$ be sets of integers of the same size, and let $M$ be a matrix. Then $M_{U,W}$ is the submatrix of $M$ with the rows corresponding to $U$ removed and the columns corresponding to $W$ removed.
\end{definition}

\begin{theorem}
\textbf{Lewis Carroll Identity.} Let $M$ be a square matrix. Then
$$\det(M)\cdot\det(M_{\{1,2\},\{1,2\}})=\det(M_{\{2\},\{2\}})\cdot\det(M_{\{1\},\{1\}})-\det(M_{\{2\},\{1\}})\cdot\det(M_{\{1\},\{2\}}).$$
\end{theorem}

As a shorthand, we sometimes drop the set notation in the subscripts.

\begin{example}
Let $M=\bordermatrix{&1'&2'&3'&4'\cr 1&3&7&0&0\cr 2&8&1&0&0\cr 3&0&0&4&0\cr 4&0&0&0&2}$. Then $\det(M)=-424$, $\det(M_{12,12})=8$, $\det(M_{2,2})=24$, $\det(M_{1,1})=8$, $\det(M_{2,1})=56$, and $\det(M_{1,2})=64$, so our identity gives us $-424\cdot 8=24\cdot 8-56\cdot 64$.
\end{example}

We will interpret the Lewis Carroll Identity in terms of pairs of directed rooted forests, called the Forest Identity. Let us first define what we mean by directed rooted forests.

\begin{definition}
A \textbf{tree} is a connected directed graph with no directed cycles such
that no node has more than one edge coming out of it. A \textbf{root} is a node in a tree
such that there exists a directed path from all other nodes to the root. Note that a
root cannot have any out-edges, or else it would form a directed cycle. A \textbf{forest} is a
disjoint union of trees.
\end{definition}

Notice that, in using these definitions, a forest is a directed graph with no cycles
in which every node has either one out-edge or no out-edges. The roots in the forest
are exactly those nodes with no out-edges.

\begin{definition}
A path from node $i$ to node $j$ is called a \textbf{meta-edge} $i\to j$
\end{definition}

\begin{example}
In the following example, the forest on the left is a two-forest (forest with two trees), with one tree rooted at node $0$ and one tree rooted at node $2$. The edges $1\to 4$, $4\to 3$, and $3\to 2$ together form the meta-edge $1\to 2$. The forest on the right is also a two-forest, with one tree rooted at node $0$ and one tree rooted at node $1$. The edge $2\to 1$ forms a meta-edge $2\to 1$. Together these forests form a pair of directed rooted forests.

\begin{center}
    \begin{tikzpicture}[baseline=0ex]
\node (0) at(0,0) {0};
\node (1) at(1,1) {1};
\node (2) at(2,0) {2};
\node (3) at(1.6,-1) {3};
\node (4) at(.4,-1) {4};
\draw[thick, red, dashed, ->] (4) edge (3);
\draw[thick, red, dashed, ->] (3) edge (2);
\draw[thick, red, dashed, ->] (1) edge (4);

\node (0) at(3,0) {0};
\node (1) at(4,1) {1};
\node (2) at(5,0) {2};
\node (3) at(4.6,-1) {3};
\node (4) at(3.4,-1) {4};
\draw[thick, ->] (4) edge (0);
\draw[thick, ->] (3) edge (0);
\draw[thick,red,dashed, ->] (2) edge (1);

\end{tikzpicture}
\end{center}

\end{example}

\begin{definition}
The \textbf{weight} of an edge is a formal variable associated with the
edge. Often, we can think of the weight as a nonnegative real number. In this case,
if there is no edge $i\to j$, then we say that the weight of $i\to j$ is zero.
\end{definition}

\begin{definition}
Given a forest F, the \textbf{forest monomial} $a_F$ is the product of the
variables $a_{ij}$ representing edge weights of the edges $i\to j$ in F.
\end{definition}

\begin{example}
Let $F$ be the forest below:

\begin{center}
   \begin{tikzpicture}[baseline=0ex]
\node (0) at(0,0) {0};
\node (1) at(1,1) {1};
\node (2) at(2,0) {2};
\draw[thick, ->] (2) edge (0);
\draw[thick, ->] (1) edge (0);
\end{tikzpicture}
\end{center}
The edge weights appearing in this forest are $a_{10}$ and $a_{20}$, so the forest monomial $a_F$ is $a_{10}a_{20}$.
\end{example}

To give the Forest Identity, we will first introduce some notation. For a set of nodes $U$, let $R_U$ be the set of forests where the nodes in $U$ (and only those nodes) are the roots. Thus $R_0$ is the set of all trees with root $0$. For an added superscript $i\to j$, let $R_U^{i\to j}$ be the set of forests with roots at $U$ and with node $i$ in the tree rooted at $j$ (that is to say there exists a meta-edge $i\to j$). Thus $R_{0,2}^{1\to 2}$ is the set of two-forests (forests containing two trees) with roots at nodes 0 and 2 and a meta-edge $1\to 2$. We will let $R^{NF}$ stand for non-forbidden forests, which are the two-forests that are allowed as part of the Forest Identity, as defined below.

\begin{theorem}\label{forestidentity}
\textbf{Forest Identity.} Let $R^{NF}$ be the set $R_{0,2}\times R_{0,1}\setminus R_{0,2}^{1\to 2}\times R_{0,1}^{2\to 1}$. Then
\begin{equation}
\sum_{(F,G)\in R_0\times R_{0,1,2}}a_Fa_G=\sum_{(F,G)\in R^{NF}}a_Fa_G.
\end{equation}
\end{theorem}

We can also represent this diagramatically. Let a star represent a root and edges represent meta-edges. Then the following diagram represents two-forests on labeled vertices $0,1,\dots,n$ with roots at nodes $0$ and $2$ and a meta-edge $1\to 0$ (that could involve other nodes, but doesn't have to).

\begin{center}
\includegraphics[scale=.5]{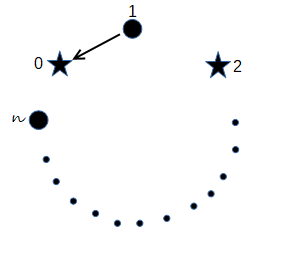}
\end{center}

Then the Forest Identity deals with pairs of forests as represented below.

\begin{center}
\includegraphics[scale=.4]{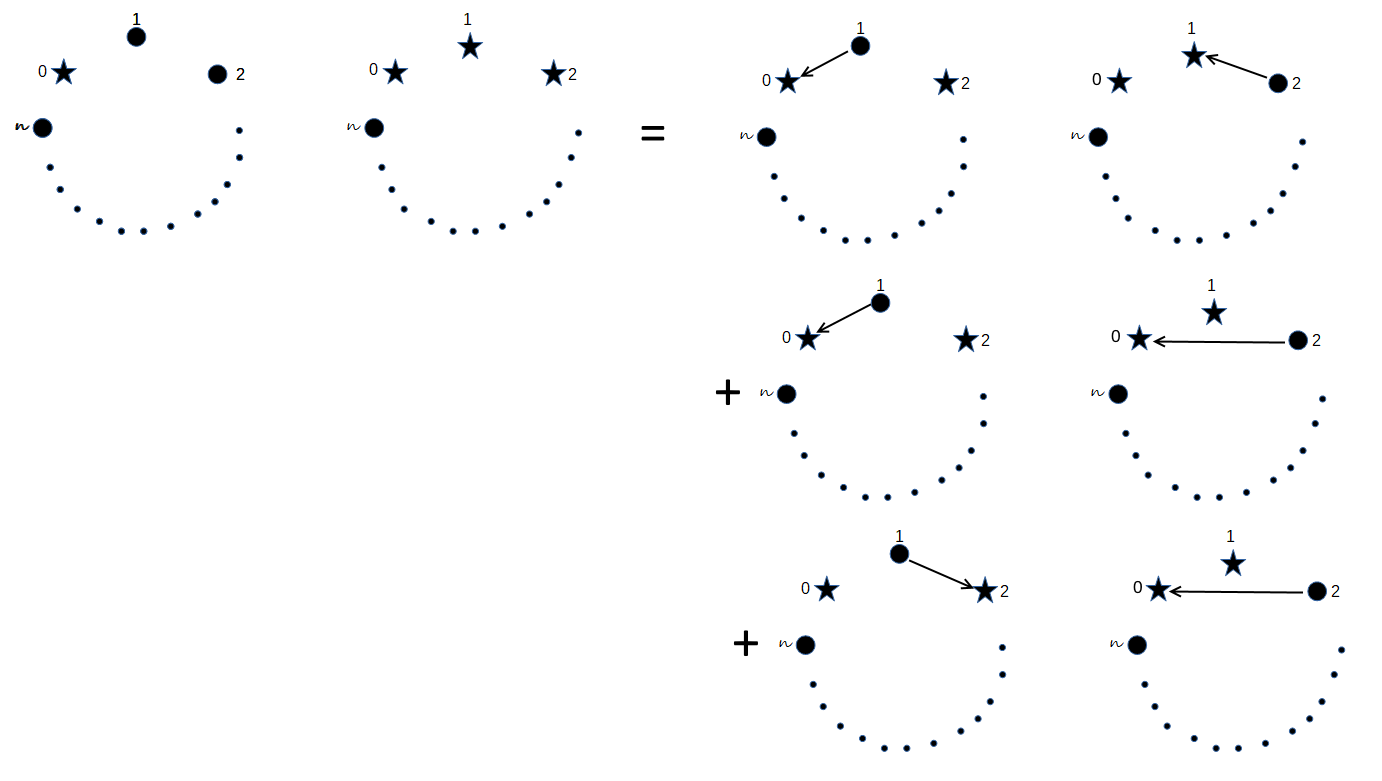}
\end{center}

\begin{example}
Below is the case $n=2$:

\begin{tikzpicture}[scale=.8][baseline=0ex]
\node (0) at(0,0) {0};
\node (1) at(1,1) {1};
\node (2) at(2,0) {2};
\draw[thick, ->] (2) edge (1);
\draw[thick, ->] (1) edge (0);

\node (0) at(3,0) {0};
\node (1) at(4,1) {1};
\node (2) at(5,0) {2};

\node (0) at(12,0) {0};
\node (1) at(13,1) {1};
\node (2) at(14,0) {2};
\draw[thick, ->] (1) edge (0);

\node (0) at(15,0) {0};
\node (1) at(16,1) {1};
\node (2) at(17,0) {2};
\draw[thick, ->] (2) edge (1);
\end{tikzpicture}

\vspace{10mm}

\begin{tikzpicture}[scale=.8][baseline=0ex]
\node (0) at(0,0) {0};
\node (1) at(1,1) {1};
\node (2) at(2,0) {2};
\draw[thick, ->] (1) edge (0);
\draw[thick, ->] (2) edge (0);

\node (0) at(3,0) {0};
\node (1) at(4,1) {1};
\node (2) at(5,0) {2};

\draw[thick, <->] (7,0) -- (10,0);

\node (0) at(12,0) {0};
\node (1) at(13,1) {1};
\node (2) at(14,0) {2};
\draw[thick, ->] (1) edge (0);

\node (0) at(15,0) {0};
\node (1) at(16,1) {1};
\node (2) at(17,0) {2};
\draw[thick, ->] (2) edge (0);
\end{tikzpicture}

\vspace{10mm}

\begin{tikzpicture}[scale=.8][baseline=0ex]
\node (0) at(0,0) {0};
\node (1) at(1,1) {1};
\node (2) at(2,0) {2};
\draw[thick, ->] (2) edge (0);
\draw[thick, ->] (1) edge (2);

\node (0) at(3,0) {0};
\node (1) at(4,1) {1};
\node (2) at(5,0) {2};

\node (0) at(12,0) {0};
\node (1) at(13,1) {1};
\node (2) at(14,0) {2};
\draw[thick, ->] (1) edge (2);

\node (0) at(15,0) {0};
\node (1) at(16,1) {1};
\node (2) at(17,0) {2};
\draw[thick, ->] (2) edge (0);
\end{tikzpicture}

\end{example}

\section[Connecting Identities]{Matrix Tree Theorem connection between Forest and Lewis Carroll Identities}

We can use the All Minors Matrix Tree Theorem \cite{chaiken} to derive the Forest Identity from the Lewis Carroll Identity.

\begin{definition}
Let $a_ij$ be the weight of the edge $i\to j$. Define the \textbf{Laplacian} $A$ by 
\[A_{ij}= \begin{cases} 
      -a_{ij} & i\neq j \\
      \displaystyle\sum_{k\neq i}a_{ik} & i=j
   \end{cases}
\]
\end{definition}

Replace $M$ in the Lewis Carroll Identity with $A_{0,0}$, where $A$ is the Laplacian matrix corresponding to the complete directed graph on nodes $\{0,1,\dots,n\}$ with a directed edge $i\to j$ for all $i\neq j$, and $A_{0,0}$ is $A$ with the zeroth row and column removed. Then we obtain the following:

\begin{equation}\label{dodgson}
\det(A_{0,0})\cdot\det(A_{012,012})=\det(A_{01,01})\cdot\det(A_{02,02})-\det(A_{02,01})\cdot\det(A_{01,02}).
\end{equation}


Now we can use the Matrix Tree Theorem \cite{chaiken} to see that when $U=W$ in $\det(A_{U,W})$, we have the sum of the weights of the forests rooted at the nodes in set $U$. When $U\neq W$, we have that $\det(A_{02,01})\cdot\det(A_{01,02})$ is pairs of forests of the form 
\begin{center}
\includegraphics[scale=.5]{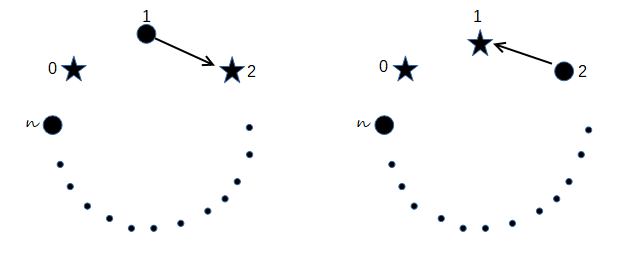}
\end{center}

Since we are subtracting these from the total number of pairs of forests rooted at $01$ and $02$, we will call these \textbf{forbidden forests}. The meta-edges $1\to 2$ and $2\to 1$ that cause a forbidden forest will together be called a \textbf{forbidden meta-cycle}. Recall from our definition of the Forest Identity that $R^{NF}$ represents the non-forbidden forests since it removes the set of forbidden forests $R^{1\to 2}_{0,2}\times R^{2\to 1}_{0,1}$ from the total set of pairs two-forests rooted at $02$ and $01$.

Taken together, we have that the left hand side of the Forest Identity corresponds to the left hand side of Eq. \ref{dodgson}, and the right hand side of the Forest Identity corresponds to the right hand side of Eq. \ref{dodgson}. The Lewis Carroll Identity thus proves the Forest Identity. Our goal for the rest of the paper is to prove the Forest Identity directly, and then use it to prove the Lewis Carroll Identity.

\section{The Red Hot Potato algorithm}

Our algorithm is based on a consequence of the Involution Principle \cite{garsia}, \cite{sally}. Before we state it, we need a few definitions:

\begin{definition}
A \textbf{signed set} is a set $S$ that is partitioned into two pieces, $S^+$ and $S^-$, such that $S = S^+\sqcup S^-$. A \textbf{sign-reversing function} on $S$ is a function that sends elements from $S^+$ to $S^-$ and elements from $S^-$ to $S^+$.
\end{definition}

\begin{definition}\label{difference}
Given two signed sets, $S_1$ and $S_2$, the \textbf{difference} $S_1 - S_2$ is the
disjoint union of the two sets such that $(S_1-S_2)^+ = S_1^+ \sqcup S^-_2$
and $(S_1-S_2)^- = S^-_1 \sqcup S^+_2$.
\end{definition}

\begin{theorem}\label{invothm}
Given any sequence of signed sets $S_0, S_1, \dots , S_{k+1}$ where $S_0$ and 
$S_{k+1}$ contain only
positive elements, and sign-reversing involutions $\phi_0,\dots, \phi_k$ where $\phi_i: S_i-S_{i+1} \to S_i-S_{i+1}$, there
is a constructible bijection between $S_0$ and $S_{k+1}$.
\end{theorem}

The Involution Principle works by applying $\phi_0$ to an element in $S_0$, and then
continuing to apply the appropriate involution to the output of each prior involution. Because every set $S_i$ is involved in the domain (and range) of two involutions ($\phi_i$ and $\phi_{i-1}$), the appropriate involution is the one of those two possibilities that was not
just used.

Our goal will be to find a sequence of signed sets and sign-reversing involutions satisfying the assumptions of Theorem \ref{invothm} such that $S_0$ is the set of pairs of trees and three-forests (the left hand side of Theorem \ref{forestidentity}), and $S_{k+1}$ is the set of non-forbidden two-forests (the right hand side of Theorem \ref{forestidentity}). In this case, $k=2$, so our set of non-forbidden two-forests will be labeled $S_3$.

\subsection{Sets}

Our signed sets will each involve pairs of graphs with edges colored either black or red (dashed).

Let $S_0$ be pairs of trees and three-forests, the trees with roots at node 0 and the three-forests with roots at nodes 0, 1, and 2 (the pairs involved in the left hand side of the Forest Identity). In this set, all edges will be black and all elements will be positive.

Let $S_1$ be pairs of graphs, one with no edge out of node $0$ and one edge out of each of the rest of the nodes, and one with no edge out of nodes $0$, $1$, or $2$ and one edge out of each of the rest of the nodes. Each cycle in either graph can be individually colored either red or black, and edges not involved in cycles must be colored black. If an even number of cycles in the two graphs is red, then the pair is positive. If an odd number of cycles is red, the pair is negative.

Let $S_2$ be pairs of graphs, one with no edge out of nodes $0$ and $2$ and one edge out of each of the rest of the nodes, and one with no edge out of nodes $0$ and $1$ and one edge out of each of the rest of the nodes. If there is a forbidden meta-cycle in the pair, we count it as a single cycle. Each cycle can either be colored red or black. The rest of the edges must be colored black. If an even number of cycles is colored red, then the pair is positive. If an odd number of cycles is colored red, the pair is negative.

Let $S_3$ be pairs of two-forests, one with roots at nodes 0 and 2 and the other with roots at nodes 0 and 1, that are not forbidden (the pairs involved in the right hand side of the Forest Identity). In this set, all edges will be black and all elements will be positive.

\subsection{Involutions}

We define involutions $\phi_0$, $\phi_1$, and $\phi_2$ on differences of sets and prove they are sign-reversing. For clarity, we will call the graph with an edge coming out of node 1 $A$ and we will call the other $B$. We will call the graph with an edge coming out of node $2$ $F$ and will call the other $R$ (in $\phi_1$ we will be ``moving'' forward and backwards along edges. We move forwards along the edges in $F$, and we reverse along the edges in $R$). In $S_1$, $A$ and $F$ are the same graph, whereas in $S_2$, $B$ and $F$ are the same.

We define $\phi_0:S_0-S_1\to S_0-S_1$ and $\phi_2:S_2-S_3\to S_2-S_3$ in the same way. Below, we will define $\phi_0$. We define $\phi_2$ in a similar way by replacing $S_0$ with $S_3$ and $S_1$ with $S_2$. Recall from Definition \ref{difference} that when we talk about $S_0-S_1$, we mean the disjoint union of $S_0$ and $S_1$ where the signs in $S_1$ have been reversed.
\begin{itemize}
\item Notice that $S_0\subset S_1$. If $t\in S_0\cap S_1$, then $\phi_0$ sends $t\in S_0$ to itself in $S_1$ and vice versa. This is clearly an involution. It is sign-reversing since $t$ is positive in both $S_0$ and $S_1$, so it is negative in $-S_1$.

\item If $t\in S_1$ and $t\notin S_0$, then there must be at least one cycle in the graphs. Then the involution changes the color of one cycle. If there is a cycle in $A$, we change the color of the cycle in $A$ containing the largest node. If there are no cycles in $A$, we change the color of the cycle in $B$ containing the largest node. If there is a forbidden meta-cycle, we change its color if there are no other cycles in the two graphs (this only applies to $\phi_2$). This is clearly an involution. It is also sign-reversing since changing the color of one cycle changes the parity of the number of red cycles.

\end{itemize}

\begin{figure}
\begin{center}
\includegraphics[scale=.45]{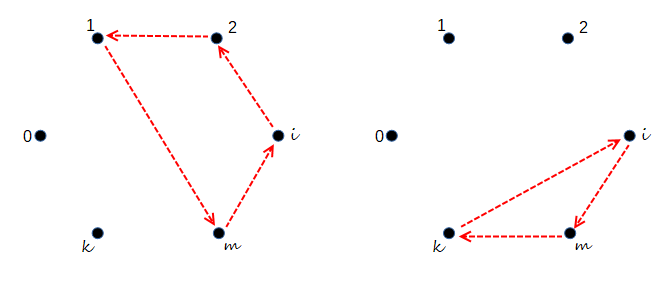}\\
\hspace{20mm}\includegraphics[scale=.25]{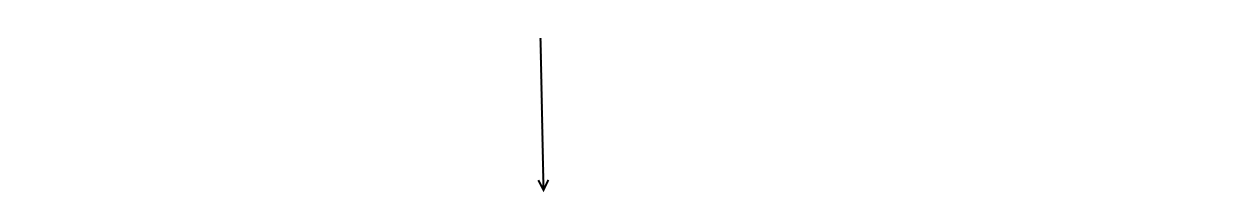}\\ 
\hspace{10mm}\includegraphics[scale=.45]{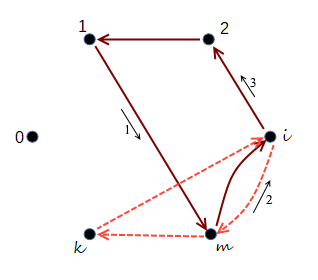}\\
\hspace{25mm}\includegraphics[scale=.25]{darrow.png}\\
\hspace{10mm}\includegraphics[scale=.45]{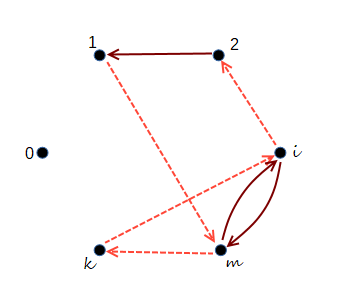}\\
\hspace{25mm}\includegraphics[scale=.25]{darrow.png}\\
\hspace{5mm}\includegraphics[scale=.45]{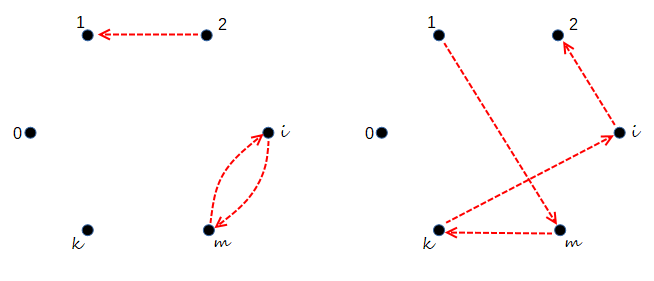}\\
\caption{An illustration of the crabwalk. All of the edges shown here are red meta-edges. The smaller arrows along the edges indicate the order in which each meta-edge is moved from one graph to another.\label{crabwalk}}
\end{center}
\end{figure}

The involution $\phi_1:S_1-S_2\to S_1-S_2$ is the involution that actually moves edges back and forth between the graphs. This part of the algorithm is based on Zeilberger's bijection, \cite{zeilberger}. When we talk about moving a red meta-edge $i\to a_1\to \dots\to a_k\to j$ from $A$ to $B$, we mean that we move all of the red edges in the $i\to j$ meta-edge from $A$ to $B$, and move all of the black edges in $B$ coming out of the nodes $a_1,a_2,\dots, a_k$ to $A$. Then $\phi_1(\{A,B\})=\{C,D\}$ where $C$ and $D$ are defined as follows:

\begin{itemize}
\item If the edge coming out of node 1 is black, move that edge from $A$ to $B$ to form a new pair $\{C,D\}$ (where now $D$ has the edge out of node 1).

\item If the edge coming out of node 1 is red, then we move the edges designated by the \textbf{crabwalk}, defined as follows. Create a graph with the same node set as $F$, and with edge set the set of all red edges from $F$ and $R$. Color the edges coming from $F$ dark red and the edges coming from $R$ light red (dashed). This is the crabwalk colored graph. We will always move forward along dark red edges (edges from $F$), and backward along light red edges (edges from $R$). 

If the edge out of node $1$ is in $F$, then we begin by moving forward along the dark red edge coming out of $1$, changing its color to light red starting with the tail of the edge and then coloring the head. Changing the color half an edge at a time becomes relevant in the proof of Lemma \ref{Parity}. We continue along that meta-edge until we reach a node that has a light red edge going into it. We travel backwards along the light red meta-edge, changing first the head of the edge to dark red and then the tail to dark red, until we reach a node that has a dark red edge coming out of it. Then we travel forward along the dark red meta-edge, changing it to light red, until we reach a node that has a light red edge going into it. We continue in this manner until we have reached either node 2 or node 1.

If the edge out of node $1$ is in $R$, then let $i\in \{1,2\}$ be the node that $1$ is rooted at, that is to say that there is a meta-edge $1\to i$. Then we begin by moving backwards along the light red edge coming into node $i$, changing first the head of the edge to dark red and then the tail. We continue moving backwards along the meta-edge until we reach a node that has a dark red edge coming out of it. Then we travel forward along the dark red meta-edge, changing it to light red, and so on until we have reached either node 2 or node 1.

Returning to graphs $F$ and $R$, we move the red meta-edges that changed color in the crabwalk, so that all dark red edges are in $F$ and all light red ones are in $R$. Figure \ref{crabwalk} gives an example of the crabwalk.

\end{itemize}

We will prove that $\phi_1$ is a sign-reversing involution in Section 5.
Assuming that $\phi_1$ is in fact a sign-reversing involution, then these sets and involutions satisfy the hypotheses of Theorem \ref{invothm}, so we have proved the Forest Identity. Following the algorithm for finding the bijection that the Involution Principle guarantees (see \cite{sally}), we can construct the bijection as follows:

\textbf{Red Hot Potato algorithm:}

We begin with a pair of forests in $S_0$. We apply $\phi_0$ and then apply $\phi_1$. Once finished with $\phi_1$, we change the color of the appropriate cycle using $\phi_2$ and do $\phi_1$ again. Then we change the color of the appropriate cycle and so on. We finish when, upon performing $\phi_2$, there are no colors to be changed, i.e. when $\phi_2$ yields a pair of graphs in $S_3$.

\section{Example}

\begin{tikzpicture}[scale=.9][baseline=0ex]
\node (0) at(0,0) {0};
\node (1) at(1,1) {1};
\node (2) at(2,0) {2};
\node (3) at(1.6,-1) {3};
\node (4) at(.4,-1) {4};
\draw[thick, ->] (4) edge (3);
\draw[thick, ->] (3) edge (0);
\draw[thick, ->] (2) edge (1);
\draw[thick, ->] (1) edge (4);

\node (0) at(3,0) {0};
\node (1) at(4,1) {1};
\node (2) at(5,0) {2};
\node (3) at(4.6,-1) {3};
\node (4) at(3.4,-1) {4};
\draw[thick, ->] (4) edge (2);
\draw[thick, ->] (3) edge (4);
\end{tikzpicture}\\ Applying $\phi_0$ simply returns the same pair in $S_1$.

\begin{tikzpicture}[scale=.9][baseline=0ex]
\node (0) at(0,0) {0};
\node (1) at(1,1) {1};
\node (2) at(2,0) {2};
\node (3) at(1.6,-1) {3};
\node (4) at(.4,-1) {4};
\draw[thick, ->] (4) edge (3);
\draw[thick, ->] (3) edge (0);
\draw[thick, ->] (2) edge (1);

\node (0) at(3,0) {0};
\node (1) at(4,1) {1};
\node (2) at(5,0) {2};
\node (3) at(4.6,-1) {3};
\node (4) at(3.4,-1) {4};
\draw[thick, ->] (4) edge (2);
\draw[thick, ->] (3) edge (4);
\draw[thick, ->] (1) edge (4);

\draw[thick, ->] (7,0) -- (8,0);

\node (0) at(10,0) {0};
\node (1) at(11,1) {1};
\node (2) at(12,0) {2};
\node (3) at(11.6,-1) {3};
\node (4) at(10.4,-1) {4};
\draw[thick, ->] (4) edge (3);
\draw[thick, ->] (3) edge (0);
\draw[thick, red, dashed, ->] (2) edge (1);

\node (0) at(13,0) {0};
\node (1) at(14,1) {1};
\node (2) at(15,0) {2};
\node (3) at(14.6,-1) {3};
\node (4) at(13.4,-1) {4};
\draw[thick, red, dashed, ->] (4) edge (2);
\draw[thick, ->] (3) edge (4);
\draw[thick, red, dashed, ->] (1) edge (4);
\end{tikzpicture} \\ The left is the result of $\phi_1$ in $S_2$, moving the single black edge out of $1$ from $A$ to $B$. The right is the result of $\phi_2$, which changes the color of the forbidden meta-cycle.

\begin{tikzpicture}[scale=.9][baseline=0ex]
\node (0) at(0,0) {0};
\node (1) at(1,1) {1};
\node (2) at(2,0) {2};
\node (3) at(1.6,-1) {3};
\node (4) at(.4,-1) {4};
\draw[thick, red, dashed, ->] (4) edge (2);
\draw[thick, ->] (3) edge (0);
\draw[thick, red, dashed, ->] (2) edge (1);
\draw[thick, red, dashed, ->] (1) edge (4);

\node (0) at(3,0) {0};
\node (1) at(4,1) {1};
\node (2) at(5,0) {2};
\node (3) at(4.6,-1) {3};
\node (4) at(3.4,-1) {4};
\draw[thick, ->] (4) edge[bend right] (3);
\draw[thick, ->] (3) edge[bend right] (4);

\draw[thick, ->] (7,0) -- (8,0);

\node (0) at(10,0) {0};
\node (1) at(11,1) {1};
\node (2) at(12,0) {2};
\node (3) at(11.6,-1) {3};
\node (4) at(10.4,-1) {4};
\draw[thick, ->] (4) edge (2);
\draw[thick, ->] (3) edge (0);
\draw[thick, ->] (2) edge (1);
\draw[thick, ->] (1) edge (4);

\node (0) at(13,0) {0};
\node (1) at(14,1) {1};
\node (2) at(15,0) {2};
\node (3) at(14.6,-1) {3};
\node (4) at(13.4,-1) {4};
\draw[thick, ->] (4) edge[bend right] (3);
\draw[thick, ->] (3) edge[bend right] (4);
\end{tikzpicture} \\ The left is the result of $\phi_1$, where we moved backwards along the meta-edge $1\to 2$. The right is the result of $\phi_2$, which changes the color of the cycle in $A$.

\begin{tikzpicture}[scale=.9][baseline=0ex]
\node (0) at(0,0) {0};
\node (1) at(1,1) {1};
\node (2) at(2,0) {2};
\node (3) at(1.6,-1) {3};
\node (4) at(.4,-1) {4};
\draw[thick, ->] (4) edge (2);
\draw[thick, ->] (3) edge (0);
\draw[thick, ->] (2) edge (1);

\node (0) at(3,0) {0};
\node (1) at(4,1) {1};
\node (2) at(5,0) {2};
\node (3) at(4.6,-1) {3};
\node (4) at(3.4,-1) {4};
\draw[thick, ->] (4) edge[bend right] (3);
\draw[thick, ->] (3) edge[bend right] (4);
\draw[thick, ->] (1) edge (4);

\draw[thick, ->] (7,0) -- (8,0);

\node (0) at(10,0) {0};
\node (1) at(11,1) {1};
\node (2) at(12,0) {2};
\node (3) at(11.6,-1) {3};
\node (4) at(10.4,-1) {4};
\draw[thick, ->] (4) edge (2);
\draw[thick, ->] (3) edge (0);
\draw[thick, ->] (2) edge (1);

\node (0) at(13,0) {0};
\node (1) at(14,1) {1};
\node (2) at(15,0) {2};
\node (3) at(14.6,-1) {3};
\node (4) at(13.4,-1) {4};
\draw[thick, red, dashed, ->] (4) edge[bend right] (3);
\draw[thick, red, dashed, ->] (3) edge[bend right] (4);
\draw[thick, ->] (1) edge (4);
\end{tikzpicture} \\ The left is the result of $\phi_1$ in $S_2$, moving the single black edge out of node $1$. The right is the result of $\phi_2$, which changes the color of the cycle in $A$ (which in this case is the right-hand graph).

\includegraphics[scale=.675]{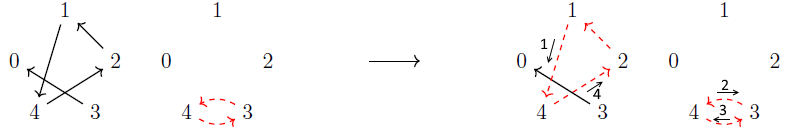}

The left is the result of $\phi_1$ in $S_1$, moving the single black edge out of node $1$. The right is the result of $\phi_2$, which changes the color of the cycle in $A$.

\begin{tikzpicture}[scale=.9][baseline=0ex]
\node (0) at(0,0) {0};
\node (1) at(1,1) {1};
\node (2) at(2,0) {2};
\node (3) at(1.6,-1) {3};
\node (4) at(.4,-1) {4};
\draw[thick, red, dashed, ->] (2) edge (1);
\draw[thick, red, dashed, ->] (4) edge[bend right] (3);
\draw[thick, red, dashed, ->] (3) edge[bend right] (4);

\node (0) at(3,0) {0};
\node (1) at(4,1) {1};
\node (2) at(5,0) {2};
\node (3) at(4.6,-1) {3};
\node (4) at(3.4,-1) {4};
\draw[thick, red, dashed, ->] (1) edge (4);
\draw[thick, red, dashed, ->] (4) edge (2);
\draw[thick, ->] (3) edge (0);

\draw[thick, ->] (7,0) -- (8,0);

\node (0) at(10,0) {0};
\node (1) at(11,1) {1};
\node (2) at(12,0) {2};
\node (3) at(11.6,-1) {3};
\node (4) at(10.4,-1) {4};
\draw[thick, red, dashed, ->] (2) edge (1);
\draw[thick, ->] (4) edge[bend right] (3);
\draw[thick, ->] (3) edge[bend right] (4);

\node (0) at(13,0) {0};
\node (1) at(14,1) {1};
\node (2) at(15,0) {2};
\node (3) at(14.6,-1) {3};
\node (4) at(13.4,-1) {4};
\draw[thick, red, dashed, ->] (1) edge (4);
\draw[thick, red, dashed, ->] (4) edge (2);
\draw[thick, ->] (3) edge (0);
\end{tikzpicture} \\ The left is the result of $\phi_1$ in $S_2$. We began in $F$, moving the edge $1\to 4$. This intersected a red cycle in $R$, so we then moved $3\to 4$ followed by $4\to 3$ (since we move backwards in $R$). Upon returning to $4$, we intersected the red cycle in $F$, so we moved $4\to 2$ and, on reaching $2$, finished the crabwalk. The right is the result of $\phi_2$, which changes the color of the cycle in $B$.

\begin{tikzpicture}[scale=.9][baseline=0ex]
\node (0) at(0,0) {0};
\node (1) at(1,1) {1};
\node (2) at(2,0) {2};
\node (3) at(1.6,-1) {3};
\node (4) at(.4,-1) {4};
\draw[thick, red, dashed, ->] (1) edge (4);
\draw[thick, red, dashed, ->] (4) edge (2);
\draw[thick, red, dashed, ->] (2) edge (1);
\draw[thick, ->] (3) edge (4);

\node (0) at(3,0) {0};
\node (1) at(4,1) {1};
\node (2) at(5,0) {2};
\node (3) at(4.6,-1) {3};
\node (4) at(3.4,-1) {4};
\draw[thick, ->] (3) edge (0);
\draw[thick, ->] (4) edge (3);

\draw[thick, ->] (7,0) -- (8,0);

\node (0) at(10,0) {0};
\node (1) at(11,1) {1};
\node (2) at(12,0) {2};
\node (3) at(11.6,-1) {3};
\node (4) at(10.4,-1) {4};
\draw[thick, ->] (1) edge (4);
\draw[thick, ->] (4) edge (2);
\draw[thick, ->] (2) edge (1);
\draw[thick, ->] (3) edge (4);

\node (0) at(13,0) {0};
\node (1) at(14,1) {1};
\node (2) at(15,0) {2};
\node (3) at(14.6,-1) {3};
\node (4) at(13.4,-1) {4};
\draw[thick, ->] (3) edge (0);
\draw[thick, ->] (4) edge (3);
\end{tikzpicture} \\ The left is the result of $\phi_1$ in $S_1$, moving the meta-edge $1\to 2$. The right is the result of $\phi_2$, which changes the color of the cycle in $A$.

\begin{tikzpicture}[scale=.9][baseline=0ex]
\node (0) at(0,0) {0};
\node (1) at(1,1) {1};
\node (2) at(2,0) {2};
\node (3) at(1.6,-1) {3};
\node (4) at(.4,-1) {4};
\draw[thick, ->] (4) edge (2);
\draw[thick, ->] (2) edge (1);
\draw[thick, ->] (3) edge (4);

\node (0) at(3,0) {0};
\node (1) at(4,1) {1};
\node (2) at(5,0) {2};
\node (3) at(4.6,-1) {3};
\node (4) at(3.4,-1) {4};
\draw[thick, ->] (3) edge (0);
\draw[thick, ->] (4) edge (3);
\draw[thick, ->] (1) edge (4);
\end{tikzpicture}\\ The result of $\phi_1$ in $S_2$, moving the single black edge out of node $1$. When we apply $\phi_2$ to this graph, no colors change, so we end up in $S_3$. Thus we are done!

\section[Proof of $\phi_1$]{Proof that $\phi_1$ is a sign-reversing involution}

\begin{lemma}\label{inout}
For any pair of graphs in the image of $\phi_1$, if node $i$ has a red edge pointing into it for any $i\in\{3,\dots,n\}$, it will also have a red edge pointing out of it in the same graph. If node $1$ or $2$ has a red edge pointing into it, it will also have a red edge pointing out of it, though not necessarily in the same graph.
\end{lemma}

\begin{proof}
Notice that this statement is equivalent to saying that red edges are only involved in cycles (for nodes $1$ and $2$, if the red edges pointing in and out of the node are in the same graph, this makes a cycle, and if they are in different graphs, this makes a meta-cycle). Because we are starting with $\{F,R\}\in S_1-S_2$, this is true for the original pair of graphs before we apply $\phi_1$. Let $i$ be a node in $\{3,\dots,n\}$. If node $i$ has one red edge in and one out of it in exactly one of the graphs, then $i$ is not a node at which we switch colors in the crabwalk. Thus if either the red edge into or out of node $i$ gets moved, the other one must as well since they are both part of the same meta-edge, and we end with the red edge in and out of $i$ still in the same graph. 

Now suppose that node $i$ has both edges out as red (and hence two edges in as red) amongst the two graphs. Then each graph has one red edge into node $i$ and one red edge out of node $i$. We want to show that this is still true after applying $\phi_1$. Let us begin by showing that we end with one red edge out of node $i$ in each graph. If we move neither red edge out of $i$ during the crabwalk, then we are done. Suppose we move the red edge out of $i$ from $R$ to $F$. Then we intersect with a red cycle in $F$, so the next step in the crabwalk is to move the red edge out of $i$ from $F$ to $R$, and thus both edges have switched graphs. This is the only way to have moved edges out of $i$ when both of the edges are red. Thus we must end with one red edge out of node $i$ in each graph. 

Now let us show that we end with one red edge into node $i$ in each graph. If neither of the red edges into node $i$ gets moved, we are done. If a red edge into node $i$ gets moved from $F$ to $R$, then the next step of the crabwalk is to move the other red edge into $i$ from $R$ to $F$, and we have finished with both graphs having a red edge into node $i$. This is the only way to have moved edges into $i$ when both of the edges are red. Thus both graphs finish with a red edge into node $i$ and a red edge out of node $i$.

If $i$ is either $1$ or $2$, then if there is a red edge out of $i$ in the original $\{F,R\}$ before applying $\phi_1$, there must also be a red edge into $i$ (the only difference in this case is that these edges do not necessarily need to be in the same graph). Because we do not add or delete edges, we will also end with a red edge out of $i$ and into $i$.
\end{proof}

\begin{lemma}\label{Range}
Let $\{A,B\}\in S_1-S_2$. Then $\phi_1(\{A,B\})\in S_1-S_2$.
\end{lemma}

\begin{proof}
To prove that $\phi_1(\{A,B\})=\{C,D\}\in S_1-S_2$, we need to show that $\{C,D\}$ has the following two defining characteristics:
\begin{enumerate}
\item There are two edges out of each node, one in each graph, except for nodes 0, 1, and 2. Node 0 has no out-edges in either graph, and nodes 1 and 2 each have a total of one out-edge between the two graphs (if the out-edges for nodes 1 and 2 are in the same graph, we are in $S_1$, and if they are in different graphs, we are in $S_2$).
\item Red edges are only in cycles (this includes the forbidden meta-cycle).
\end{enumerate}

Requirement 2 is already proved by Lemma \ref{inout}. Thus we need only show requirement 1. Because $\phi_1$ simply moves edges around, and does not add or delete any edges, then since we have started with no edges out of node 0, one edge out of nodes 1 and 2, and two edges out of the rest, we will end with that as well. Thus we just need to show that for the two edges out of nodes $3,\dots, n$ in $\{C,D\}$, one is in each graph. Since we begin with one edge out of a node $i$ in each graph, we want to show that, if we move one edge out of $i$ from $F$ to $R$, then we must move the other edge out of $i$ from $R$ to $F$. Notice that if both edges out of node $i$ are black, then neither will get moved, so they end as they started with one in each graph. If one edge is red and the other is black then, by how we defined moving red meta-edges, if the red edge gets moved, then the black one will as well. Finally, suppose that both edges out of $i$ are red. As we saw in the proof of Lemma \ref{inout}, we will end with one red edge out of $i$ in each graph.

\end{proof}

\begin{lemma}\label{path}
Every meta-edge moved from $F$ to $R$ in the course of the crabwalk joins a path in $R$ starting at $1$.
\end{lemma}

\begin{proof}
Let us break this proof into two cases:

\vspace{5mm}

\textbf{Case 1:} Let the original pair of graphs $\{F,R\}$ be in $S_1$. Then we begin by moving forward along the red meta-edge out of $1$, until we reach the first intersection point $a_1$ where $R$ had a red meta-edge pointing into it. That red meta-edge belongs to a red cycle in $R$, so when we move backwards along that cycle to the next intersection point $a_2$, we leave in $R$ the meta-edge $a_1\to a_2$ to join with the recently moved meta-edge $1\to a_1$. Thus we now have a path $1\to a_2$ in $R$. When we reach the intersection point $a_2$, our next step is to move forward along a meta-edge $a_2\to a_3$ from $F$ to $R$, so that meta-edge becomes connected to our path out of $1$. We can continue in this way showing that each meta-edge moved from $F$ to $R$ gets connected to the path from $1$.

\vspace{5mm}

\textbf{Case 2:} Let the original pair of graphs $\{F,R\}$ be in $S_2$. Then we begin by moving backwards along the edge $a_1\to i$ where $i$ is the end of the meta-edge $1\to i$ and $a_1$ is the first intersection point. Then we move forwards along the meta-edge $a_1\to a_2$, moving this edge from $F$ to $R$. But now $a_1\to a_2$ joins $1\to a_1\in R$, creating the path $1\to a_2$. We can now use the same reasoning as Case 1 to show that every meta-edge moved from $F$ to $R$ will join the path in $R$ starting at $1$.
\end{proof}

\begin{lemma}\label{involution}
The function $\phi_1:S_1-S_2\to S_1-S_2$ is an involution.
\end{lemma}

\begin{proof}
If the edge coming out of $1$ in $A$ is black, then we simply move that edge over, which means if we apply $\phi_1$ again, we just move that edge back again. Therefore in this case, $\phi_1$ is an involution.

If the edge coming out of $1$ in $A$ is red, then we argue that if we apply $\phi_1$ twice, each time the sequence of edges traveled along in the crabwalk will be the same. Since we are simply changing the colors of the edges in that sequence of edges, it is clear this would mean that $\phi_1$ is an involution.

\vspace{5mm}

\textbf{Case 1:}
Suppose that our original pair of graphs $\{F,R\}$ is in $S_1$. Since $R$ has no edges out of nodes $1$ or $2$, there can be no red edges in $R$ involving either of these nodes, so the crabwalk must end in $F$ with an edge $e_1$ going either into node $1$ or into node $2$. Thus we have moved over the edge out of $1$, but not the edge out of $2$, so $\{C,D\}\in S_2$. When applying $\phi_1$ to $\{C,D\}\in S_2$, we want to show that we begin by moving backwards along the same edge $e_1$ that we moved forward along to finish the iteration of $\phi_1$ that resulted in $\{C,D\}$. If this is the case, then by construction we will move backwards along the same sequence of edges that we originally moved forwards along, and so the same set of edges will get moved and $\phi_1(\{C,D\})$ will return to us our original $\{F,R\}$. Thus we need to show that the meta-edge $1\to i\in D$ ends with edge $e_1$. By Lemma \ref{path}, every edge moved from $F$ to $R$ becomes part of the path out of $1$. Since we finish with $e_1$ getting moved from $F$ to $R$, then $e_1$ is also part of this path. Thus when we apply $\phi_1$ to $\{C,D\}$, the first thing we do is move backwards along $e_1$ from $D$ to $C$. Thus $\phi_1(\{C,D\})=\{F,R\}$ and so $\phi_1$ is an involution.

\vspace{5mm}

\textbf{Case 2:}
Now suppose that our original pair of graphs $\{F,R\}$ is in $S_2$. If we end with $\{C,D\}\in S_1$, then we must have ended by moving the edge out of $1$ so that the edges out of nodes $1$ and $2$ end up in the same graph. Thus we must have finished the crabwalk by going backwards along the edge $e_2$ out of $1$. When we apply $\phi_1$ again, since the edges out of $1$ and $2$ are in the same graph, we begin by traveling forwards along the edge out of $1$. Since we begin by traveling the opposite direction along the last edge moved, the entire sequence will be the same by construction, and so $\phi_1$ is an involution.

Finally, let us suppose that $\{F,R\}\in S_2$ and $\{C,D\}\in S_2$. Then we must have finished by going forwards along an edge $e_3$ going into $1$ or into $2$ in $F$, since going backwards along the edge out of $1$ in $R$ would result in $\{C,D\}\in S_1$, and there is no edge out of $2$ in $R$ to move backwards along. Thus we want to show that when we apply $\phi_1$ to $\{C,D\}$, we begin by moving backwards along the last edge $e_3$ that we moved from $F$ to $R$. Since $\{C,D\}\in S_2$, we already know we will start by traveling backwards, so this stipulation is met. By Lemma \ref{path}, $e_3$ is part of the path starting at $1$. Thus when we apply $\phi_1$ to $\{C,D\}$ the first edge we move along will be backwards along $e_3$. Again, by construction, since we start by moving in the opposite direction along the last edge moved, the sequence of edges will be the same, so $\phi_1$ is an involution.
\end{proof}

\begin{lemma}\label{Parity}
If the crabwalk ends in graph $A$, then the parity of cycles (whether there are an odd or even number of cycles) remains the same after applying $\phi_1$. If the crabwalk ends in graph $B$, then the parity of cycles changes after applying $\phi_1$.
\end{lemma}

\begin{figure}[H]
\begin{subfigure}{\textwidth}
\begin{center}
\includegraphics[scale=.7]{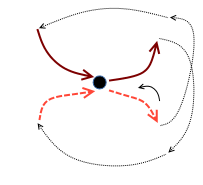}
\includegraphics[scale=.4]{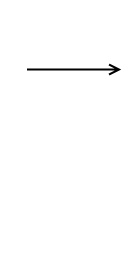}
\includegraphics[scale=.7]{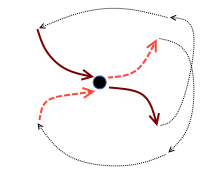}
\end{center}
\end{subfigure}\\
\begin{subfigure}{\textwidth}
\begin{center}
\includegraphics[scale=.7]{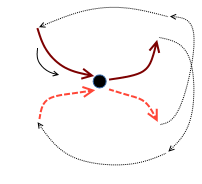}
\includegraphics[scale=.4]{arrow.png}
\includegraphics[scale=.7]{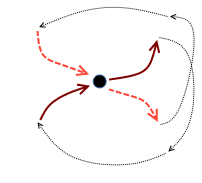}
\end{center}
\end{subfigure}
\caption{The top figure shows what happens approaching the color change along the light red edge, whereas the bottom figure shows approaching it from the dark red edge. On the left hand side, all four edges are involved in the same cycle. After the color change (on the right), that single cycle has split into two separate cycles. \label{SingleCycle}}
\end{figure}

\begin{figure}[H]
\begin{subfigure}{\textwidth}
\begin{center}
\includegraphics[scale=.7]{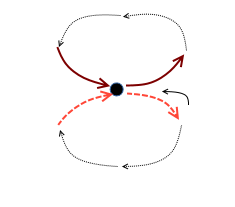}
\includegraphics[scale=.4]{arrow.png}
\includegraphics[scale=.7]{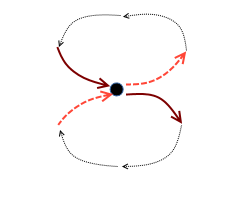}
\end{center}
\end{subfigure}\\
\begin{subfigure}{\textwidth}
\begin{center}
\includegraphics[scale=.7]{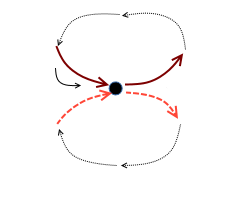}
\includegraphics[scale=.4]{arrow.png}
\includegraphics[scale=.7]{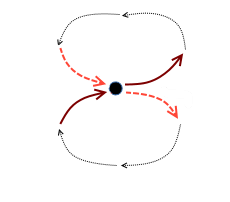}
\end{center}
\end{subfigure}
\caption{On the left hand side, there are two cycles, one involving the dark red edges and one involving the light red edges. After the color change (on the right), the two cycles have merged into a single cycle. \label{DoubleCycle}}
\end{figure}

\begin{proof}
We claim that each time we switch shades in the crabwalk (that is to say we change from moving along edges in $F$ to edges in $R$ or vice versa), the parity of the cycles changes. We always begin the crabwalk in $A$, so if we end in $A$, we have switched shades an even number of times, so the parity is the same after applying $\phi_1$. If we end in $B$, we have switched shades an odd number of times, so the parity has changed after applying $\phi_1$. 

To choose which cycle a particular edge is in partway through the crabwalk, we follow a path beginning with that edge. At every node, we follow the same shade out as we followed in, and we follow each edge its whole length regardless of shade. Because we change edge shades half an edge at a time, a cycle can consist of both shades. When we change shades, either both edges into the node have changed shades or both edges out have. Thus any node will have at most one edge in and one edge out of each shade, meaning that our method of finding cycles is well-defined. Since every edge belongs to a unique cycle, there is a well-defined number of cycles when we switch shades.

Now that we have a way of counting cycles partway through the crabwalk, let us zoom in on a node at which we will change shades. This node will have two edges (one dark red, one light red) pointing into it and pointing out of it. We will be approaching this node either forwards along the dark red edge pointing into the node, or backwards along the light red edge pointing out of the node. If we are moving forward along the dark edge into the node, then the next step in the crabwalk moves backward along the light edge into the node. Thus both edges \textbf{into} the node change shade. If we are moving backward along the light edge out of the node, then the next step is to move forward along the dark edge out of the node. Thus both edges \textbf{out of} the node change shade. Before we change shades, there are two options for the cycles that the edges are involved with.

\vspace{5mm}

\textbf{Case 1:} All four edges are involved in the same cycle. Then after the shade change, the single cycle will split into two different cycles (\textbf{Figure \ref{SingleCycle}}). Since this is the only cycle impacted by the shade change, we have changed the parity of cycles.

\vspace{5mm}

\textbf{Case 2:} The light edges are involved in one cycle and the dark edges are involved in another, so that there are two cycles before the color change. Then after the shade change, these two cycles will merge into one cycle (\textbf{Figure \ref{DoubleCycle}}). Since these are the only cycles impacted by the shade change, we have changed the parity of cycles.

Since each time we change the shade, the parity of cycles changes, then as discussed at the beginning of the proof, we have proved our lemma.
\end{proof}

\begin{theorem}
The function $\phi_1$ is a sign-reversing involution on $S_1-S_2$.
\end{theorem}

\begin{figure}[H]
\begin{center}
\includegraphics[scale=.5]{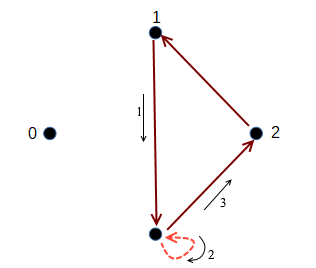}
\includegraphics[scale=.4]{arrow.png}
\includegraphics[scale=.5]{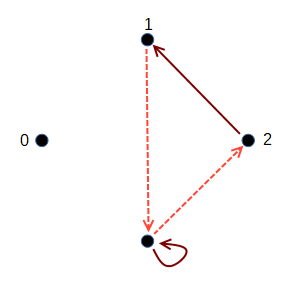}
\caption{An example of Case 1. Here, we have switched colors twice at the bottom node \label{case1}}
\end{center}
\end{figure}

\begin{figure}[H]
\begin{subfigure}{\textwidth}
\begin{center}
\includegraphics[scale=.5]{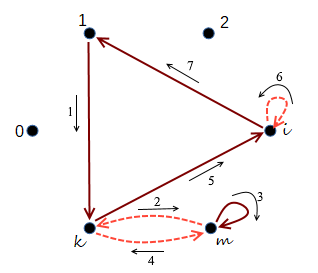}
\includegraphics[scale=.4]{arrow.png}
\includegraphics[scale=.5]{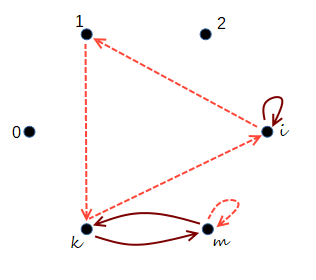}
\end{center}
\caption{An example of the crabwalk ending in $A$. The colors have switched once at $m$ and $i$ and twice at $k$, for a total of four switches. \label{case2a}}
\end{subfigure}\\
\begin{subfigure}{\textwidth}
\begin{center}
\includegraphics[scale=.5]{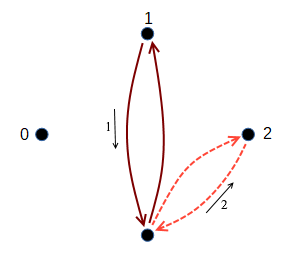}
\includegraphics[scale=.4]{arrow.png}
\includegraphics[scale=.5]{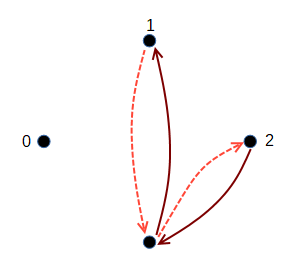}
\end{center}
\caption{An example of the crabwalk ending in $B$. The colors have switched once at the bottom node. \label{case2b}}
\end{subfigure}
\caption{Examples of Case 2.}
\end{figure}

\begin{proof}

Lemma \ref{Range} shows that $\phi_1$ is indeed a function into the correct range. Lemma \ref{involution} shows that $\phi_1$ is an involution. Thus we just need to show that $\phi_1$ is sign-reversing.

If the edge coming out of node 1 in $A$ is black, if $\{A,B\}\in S_1$, then moving the single black edge over puts $\{C,D\}\in S_2$ and vice versa. The parity of cycles does not change, which means that in $S_1-S_2$, $\phi_1$ is sign-reversing.

Now suppose that the edge coming out of node 1 in $A$ is red. Then we perform the crabwalk. We break our proof that $\phi_1$ is sign-reversing into two cases:

\vspace{5mm}

\textbf{Case 1} Let $\{A,B\}\in S_1$. Since $B$ has no edges out of nodes $1$ or $2$, there can be no red edges in $B$ involving either of these nodes, so the crabwalk must end in $A$. Then by Lemma \ref{Parity} the parity of cycles has not changed. Since we have moved over the edge out of $1$, but not the edge out of $2$, then $\{C,D\}\in S_2$, and since the parity of cycles has not changed, $\{C,D\}$ has the same sign in $S_2$ as $\{A,B\}\in S_1$. Thus the sign has changed in $S_1-S_2$. (See Figure \ref{case1} for an example).

\vspace{5mm}

\textbf{Case 2} Let $\{A,B\}\in S_2$.
\begin{itemize}
\item Suppose the crabwalk ends in $A$. Then by Lemma \ref{Parity} the parity of the cycles remains the same. Since we are not ending at node 2 in $B$, the edge out of $2$ remains in $D$, and the edge out of $1$ switches to $D$. Thus $\{C,D\}\in S_1$ and, since the sign of $\{A,B\}\in S_2$ is the same as that of $\{C,D\}\in S_1$, the sign is reversed in $S_1-S_2$. (See Figure \ref{case2a} for an example).

\item Suppose the crabwalk ends in $B$. Then by Lemma \ref{Parity} the parity of the cycles has switched. Since we move backwards along light red edges, then the last light edge moved is an edge coming out of a node. Since there is no edge out of $1$ in $B$, then we must end with the edge coming out of node $2$. Thus the edge out of $1$ has been moved to $D$ and the edge out of node $2$ has been moved to $C$, so $\{C,D\}\in S_2$. Since the parity of the cycles has been switched, the sign is reversed in $S_2$, so it is reversed in $S_1-S_2$. (See Figure \ref{case2b} for an example).
\end{itemize}

\end{proof}

\section{Asymptotic Observations}

The author programmed the Red Hot Potato algorithm into Mathematica. In this program the user can enter a pair of forests from either $S_0$ or $S_3$ and the program will generate a list of graphs similar to the list we saw in Section 4, ending with the corresponding pair of forests in $S_3$ or $S_0$ (respectively). Using this program, we observed that in more than half of cases, the algorithm took only one step to complete, namely moving the edge coming out of node 1 from graph $A$ to graph $B$. Data from running the algorithm indicates that the proportion of single step algorithms is asymptotically $2/3$. The following is a proof of this observation, kindly supplied by David Wagner from the University of Waterloo.

\begin{theorem}
The probability that the Red Hot Potato algorithm moves only a single edge is asymptotically $2/3$.
\end{theorem}

\begin{proof}
Begin with an element of $S_0$. Then $A$ is a tree rooted at $0$ and $B$ is a forest composed of three disjoint trees, which we will call $T_0$, $T_1$, and $T_2$, where the subscript indicates the root of the tree. Let $1\to v$ be the edge out of node $1$ in $A$. If $v\in T_1$ in $B$, then moving the edge $1\to v$ from $A$ to $B$ will create a cycle in $B$, requiring us to continue moving edges. On the other hand, if $v\notin T_1$, then moving the edge $1\to v$ from $A$ to $B$ will simply combine two of the trees in $B$, resulting in an element in $S_3$. Thus the Red Hot Potato algorithm will move the edge $1\to v$ (and no other edges) if and only if $v\notin T_1$. Therefore we aim to prove that $P(v\notin T_1)\to 2/3$.

Notice that $v$ cannot be 1, since that would create a loop in the original tree. If $v=0$ or if $v=2$, then $v$ cannot be in $T_1$ since the trees in $B$ are disjoint. If $v$ is any other node, then there is an equal probability of $v$ being in any of $T_1$, $T_2$ or $T_3$, so in this case, there is a $2/3$ probability that $v\notin T_1$. Given this, we can see that
\begin{align*}
    P(v\notin T_1)=&P(v\notin T_1\mid v=0)P(v=0)\\&+P(v\notin T_1\mid v=2)P(v=2)\\&+P(v\notin T_1\mid v\neq 0,2)P(v\neq 0,2)\\
    =&1\cdot P(v=0)+1\cdot P(v=2)+\frac{2}{3}P(v\neq 0,2)
\end{align*}

Since node $0$ is the root of $A$, we must calculate the probability that $v=0$ separately. In any tree on $n$ nodes, there are $n-1$ edges, giving a total degree of $2n-2$. Then the expected degree of any given vertex is $\frac{2n-2}{n}$. Since $0$ is the root, this degree is in-degree. Those $\frac{2n-2}{n}$ edges might come from any of the $n-1$ other nodes in the graph. Then any given edge comes from node $1$ with probability $\frac{1}{n-1}$, so the edge out of $1$ is $1\to 0$ with probability $\frac{2n-2}{n}\cdot \frac{1}{n-1}$. Thus $$P(v=0)=\frac{2}{n}.$$

The probability that $v\neq 0$ is therefore $1-\frac{2}{n}=\frac{n-2}{n}$. There are $n-2$ possible nodes that aren't $0$ (or $1$ since, as established above, that is impossible). Thus $$P(v=2)=\frac{n-2}{n}\cdot \frac{1}{n-2}=\frac{1}{n}.$$

The probability that $v$ is any specific node (aside from $0$ or $1$) is similarly $\frac{1}{n}$. Given that there are $n-3$ nodes that are not $0$, $1$, or $2$, we see that $$P(v\neq 0,2)=\frac{n-3}{n}.$$

Putting this all together, we get that $$P(v\notin T_1)=\frac{2}{n}+\frac{1}{n}+\frac{2}{3}\frac{n-3}{n}=\frac{2n+3}{3n}.$$ As $n\to\infty$, $P(v\notin T_1)\to\frac{2}{3}$.
\end{proof}

\section[Identity Proofs]{Proofs of Forest and Lewis Carroll Identities}

We will restate and prove our two identities using the Red Hot Potato algorithm.

\begin{theorem}
\textbf{Forest Identity.} Let $R^{NF}$ be the set $R_{0,2}\times R_{0,1}\setminus R_{0,2}^{1\to 2}\times R_{0,1}^{2\to 1}$. Then
$$\sum_{(F,G)\in R_0\times R_{0,1,2}}a_Fa_G=\sum_{(F,G)\in R^{NF}}a_Fa_G.$$
\end{theorem}

\begin{proof}
We have proven that the signed sets $S_0$, $S_1$, $S_2$, and $S_3$, and the sign-reversing involutions $\phi_0$, $\phi_1$, and $\phi_2$ satisfy the hypotheses of Theorem \ref{invothm}. Thus there exists a bijection between $S_0$ and $S_3$. Suppose that the edges in our graphs are unweighted. Then the left hand side of the identity merely counts the number of pairs $(F,G)\in R_{0}\times R_{0,1,2}$, i.e. the left hand side counts the number of elements in $S_0$. Similarly the right hand side counts the number of pairs $(F,G)\in R^{NF}$, i.e. the number of elements in $S_3$. Since there is a bijection between $S_0$ and $S_3$, then $|S_0|=|S_3|$, so the identity holds. Now suppose that the graphs are weighted. We modify our bijection between $S_0$ and $S_3$ by first taking away the weights on the edges, performing the bijection, and then putting the weights back on the respective edges. Since the edges are only moved between graphs, not changed, the weight remains the same. In this manner, we are still matching each pair of forests in $S_0$ with one in $S_3$, and the total weight $a_Fa_G$ remains the same. Thus our identity holds.
\end{proof}

\begin{theorem}\label{lci}
\textbf{Lewis Carroll Identity.} Let $M$ be a square matrix. Then
$$\det(M)\cdot\det(M_{12,12})=\det(M_{2,2})\cdot\det(M_{1,1})-\det(M_{2,1})\cdot\det(M_{1,2}).$$
\end{theorem}

\begin{proof}
The Forest Identity proves the Lewis Carroll Identity provided that $M=A_{0,0}$ where $A$ is the Laplacian for some directed graph. Notice that provided that the row sums of a matrix $A$ are zero, then $A$ is the Laplacian for a directed graph, since given any matrix $A$ with zero row sums, we can create the corresponding directed graph by giving every edge $i\to j$ the negative of the weight found in row $i$ column $j$, where $i\neq j$. Then given any square matrix $M$, we can turn it into a Laplacian for a directed graph by adding a zeroth row and column where the zeroth row can be anything that sums to zero, and the $i$th entry of the zeroth column is the opposite of the sum of the $i$th row. Thus any square matrix is of the form $A_{0,0}$ where $A$ is the Laplacian for some directed graph.
\end{proof}

\bibliographystyle{plain}

\end{document}